\documentclass[11pt]{article}
\usepackage{amsmath,amsthm,amssymb,amsfonts}

\theoremstyle{plain}
\newtheorem{theorem}{Theorem}
\newtheorem{proposition}[theorem]{Proposition}
\newtheorem{lemma}[theorem]{Lemma}
\newtheorem{corollary}[theorem]{Corollary}
\theoremstyle{definition}
\newtheorem*{remark}{Remark}
\newtheorem*{example}{Example}
\newenvironment{myproof}[1][Proof]{\noindent\textit{#1.} }{\hfill$\Box$\\}
\numberwithin{theorem}{section}\numberwithin{equation}{section}

\input{tcilatex}
\begin{document}

\title{Jet Finslerian geometry of the conformal\\
Minkowski metric}
\author{Vladimir Balan and Mircea Neagu}
\date{}
\maketitle

\begin{abstract}
The paper develops the Finsler-like geometry on the $1$-jet space for the
jet conformal Minkowski (JCM) metric, which naturally \linebreak extends the
Minkowski metric in the Chernov-Pavlov framework. To this aim there are
determined the nonlinear connection, distinguished (d-) Cartan linear
connection, d-torsions and d-curvatures. The field geometrical gravitational
and electromagnetic d-models based on the JCM metric are discussed.
\end{abstract}



\section{Introduction}

It is obvious that our genuine physical intuition distinguishes four
dimensions in a natural correspondence with the material reality.
Consequently, the four dimensionality plays a special role in almost all
modern physical theories \cite{Mikhailov,pa2,pa3}.

On the other hand, it is a well known fact that, in order to create
Relativity Theory, Einstein used Riemannian geometry instead of classical
Euclidean geometry, the first one representing a natural
ma\-the\-ma\-ti\-cal model for \textit{local isotropic space-time}. Although
the use of Riemannian geometry was indeed a genial idea, there are recent
studies of physicists that suggest a \textit{non-isotropic} perspective of
space-time. For example, in Pavlov's works \cite{Pavlov,pa2,pa3}, the
concept of inertial body mass emphasizes the necessity to study \emph{local
non-isotropic spaces}. For the study of non-isotropic physical phenomena,
Finsler geometry proves to be adequate and proficient as
ma\-the\-ma\-ti\-cal framework.

Recent studies of Russian scholars (e.g., Asanov \cite{as1}, Garas'ko \cite%
{Garasko-Book} and Pavlov \cite{Garasko-Pavlov}, \cite{Pavlov}) emphasize
the importance of Finsler geometry, which is characterized by the total
equality in rights of all non-isotropic directions. For such a reason, in
their works is underlined the important role played in theory of space-time
structure and gravitation (as well as in unified gauge field theories) by
the $m$-root metric (\cite{Shimada,Atanasiu-Neagu,Bal-Nea-Chernov}) 
\begin{equation*}
L:TM\rightarrow \mathbb{R},\mathbb{\qquad }L(x,y)=\sqrt[m]{
a_{i_1i_2...i_{m}}(x)y^{i_1}y^{i_2}...y^{i_{m}}}.
\end{equation*}
It is known that the $1$-jet fibre bundle is a basic object in the study of
classical and quantum field theories (Olver, \cite{Olver}). For these
geometrical and physical reasons, the present paper is devoted to the
construction on the $1$-jet space $J^1(\mathbb{R},M^4)$ of the Finsler-like
geometry (together with the extended gravitational and electromagnetic
geometrical models) for the \emph{jet conformal Minkowski (JCM) metric} $%
F:J^1(\mathbb{R},M^4)\rightarrow\mathbb{R},$ defined by\footnote{%
In the following we shall reduce the domain of the constructed geometric
objects in order to ensure their existence and, where this is required,
their smoothness. As well, we shall implicitly use throughout the work the
Einstein convention of summation.} {\small \ 
\begin{equation}  \label{JCM}
F(t,x,y)=e^{\sigma (x)}\cdot \sqrt{h^{11}(t)}\cdot \sqrt{%
y_1^1y_1^2+y_1^1y_1^3+ y_1^1y_1^4+y_1^2y_1^3+y_1^2y_1^4+y_1^3y_1^4},
\end{equation}
} where $\sigma (x)$ is a smooth non-constant function on $M^4$, $h^{11}(t)$
is the dual of the Riemannian metric $h_{11}(t)$ on $\mathbb{R}$ and 
\begin{equation*}
(t,x,y)=(t,x^1,x^2,x^3,x^4,y_1^1,y_1^2,y_1^3,y_1^4)
\end{equation*}
are the coordinates of the $1$-jet space $J^1(\mathbb{R},M^4)$. These
transform by the rules: 
\begin{equation}  \label{tr-rules}
{\tilde t}={\tilde t}(t),\quad {\tilde x}^{p}={\tilde x}^{p}(x^{q}),\quad {%
\tilde y}_1^{p}=\dfrac{\partial {\tilde x}^{p}}{\partial x^{q}} \dfrac{dt}{d{%
\tilde t}}\cdot y_1^{q},\qquad p,q=\overline{1,4},
\end{equation}
where $d{\tilde t}/dt\neq 0$ and rank $(\partial {\tilde x}^{p}/\partial
x^{q})=4$.

\begin{remark}
It is easy to verify that (as emphasized in the recent studies \cite{Chernov}
and \cite{Pavlov}) the geometrical object 
\begin{equation}  \label{G-11}
G_{11}(y)\overset{def}{=}y_1^1y_1^2+y_1^1y_1^3+
y_1^1y_1^4+y_1^2y_1^3+y_1^2y_1^4+y_1^3y_1^4
\end{equation}
is a quadratic form in $y=(y_1^1,y_1^2,y_1^3,y_1^4)$, whose canonical form
is the Minkowski metric. Namely, denoting $x=(x^1,x^2,x^3,x^4)$, ${\tilde x}%
=({\tilde x}^1,{\tilde x}^2,{\tilde x}^3,{\tilde x}^4)$ and $A=$%
{\footnotesize $\left(%
\begin{array}{cccc}
{1/\sqrt{6}} & -{1/\sqrt{3}} & 1 & -{1/\sqrt{6}} \\ 
{1/\sqrt{6}} & {2/\sqrt{3}} & 0 & -{1/\sqrt{6}} \\ 
{1/\sqrt{6}} & 0 & 0 & {3/\sqrt{6}} \\ 
{1/\sqrt{6}} & -{1/\sqrt{3}} & -1 & -{1/\sqrt{6}}%
\end{array}%
\right)$}, if we apply on the product manifold $\mathbb{R}\times M^4$ the
invertible linear coordinate transformation 
\begin{equation*}
t={\tilde t},\quad ^{T}x=A\text{ }\cdot \text{ }^{T}{\tilde x},
\end{equation*}
then in the induced coordinates $({\tilde t},{\tilde x},{\tilde y})$ on $J^1(%
\mathbb{R},M^4)$, we have $^{T}y=A\;\cdot\;^{T}{\tilde y}$, and the JCM
Finslerian metric (\ref{JCM}) has the particular form, 
\begin{equation*}
F({\tilde t},{\tilde x},{\tilde y})=e^{\sigma({\tilde x}\;\cdot\;^{T}A)}%
\cdot \sqrt{h^{11}({\tilde t})}\cdot\sqrt{\left( {\tilde y}_1^1\right) ^2-
\left( {\tilde y}_1^2\right) ^2-\left( {\tilde y}_1^3\right) ^2- \left({%
\tilde y}_1^4\right) ^2}.
\end{equation*}
\end{remark}

The distinguished (d-) jet framework (\cite{MA,Nc,as1,as2}), which involves
specific geometric objects as canonical nonlinear connection, Cartan
canonical linear connection, d-torsions, d-curvatures, and their related
extended gravitational and electromagnetic geometrical models produced by an
arbitrary jet Lagrangian function 
\begin{equation*}
L:J^1(\mathbb{R},M^{n})\rightarrow \mathbb{R},
\end{equation*}
was completely treated in recent works of the authors of this paper (\cite%
{BN,Neagu-Rheon}). We point out that geometrical ideas from these works are
similar, but distinct, from those promoted by Miron and Anastasiei in the
classical Lagrangian geometry on tangent bundle (\cite{MA}). Namely, the
case of the present framework was initially stated by Asanov in \cite{as2},
and further generalized in \cite{Nc} by the second author of this paper.

In the sequel, we apply the general geometrical results from \cite{BN} and 
\cite{Neagu-Rheon} to the particular jet conformal Minkowski-metric (\ref%
{JCM}).

\section{The canonical nonlinear connection of the model}

Let $(\mathbb{R},h_{11}(t))$ be a Riemannian manifold, where $\mathbb{R}$ is
the set of real numbers. The Christoffel symbol of the Riemannian metric $%
h_{11}(t)$ is 
\begin{equation*}
\varkappa _{11}^1=\frac{h^{11}}{2}\frac{dh_{11}}{dt},\quad\mbox{where}\quad
h^{11}=(h_{11})^{-1}>0.
\end{equation*}
Let also $M^4$ be a manifold of dimension four, whose local coordinates are $%
x=(x^1,x^2,x^3,x^4)$. These manifolds produce the $1$-jet space $J^1(\mathbb{%
R},M^4)$, whose local coordinates are $(t;x;y)$, where $%
y=(y_1^1,y_1^2,y_1^3,y_1^4)$.

Let's consider on $J^1(\mathbb{R},M^4)$ the JCM metric (\ref{JCM}), whose
domain of definition consists of all values $(t;x;y)$ which satisfy the
condition $G_{11}(y)>0$, where $G_{11}$ is given by (\ref{G-11}). If we use
the notation 
\begin{equation*}
S_{[1]1}=y_1^1+y_1^2+y_1^3+y_1^4,
\end{equation*}
then the following relations are true: 
\begin{equation*}
\begin{array}{l}
\medskip G_{i1}\overset{def}{=}\dfrac{\partial G_{11}}{\partial y_1^{i}}%
=S_{[1]1}-y_1^{i}, \\ 
G_{ij}\overset{def}{=}\dfrac{\partial G_{i1}}{\partial y_1^{j}}= \dfrac{%
\partial ^2G_{11}}{\partial y_1^{i}\partial y_1^{j}}=1-\delta _{ij},%
\end{array}%
\end{equation*}
where $\delta _{ij}$ is the Kronecker symbol. Obviously, the homogeneity of
degree $2$ of the "$y$-function" $G_{11}$ (which is in fact a d-tensor on $%
J^1(\mathbb{R},M^4)$) leads to the equalities: 
\begin{equation*}
G_{i1}y_1^{i}=2G_{11},\quad G_{ij}y_1^{i}y_1^{j}=2G_{11}.
\end{equation*}
By direct computation, we get

\begin{lemma}
a) The \textit{fundamental metrical d-tensor} produced by the JCM Finslerian
metric $F$ is given by the formula 
\begin{equation*}
g_{ij}(t,x,y)=\frac{h_{11}(t)}{2}\frac{\partial ^2F^2}{\partial
y_1^{i}\partial y_1^{j}},
\end{equation*}
which in our case leads to 
\begin{equation}  \label{metric}
g_{ij}(x)=\frac{e^{2\sigma (x)}}{2}\left( 1-\delta _{ij}\right),
\end{equation}
and the matrix $g=(g_{ij})$ admits the inverse $g^{-1}=(g^{jk})$, whose
entries are 
\begin{equation*}
g^{jk}(x)=\frac{2e^{-2\sigma (x)}}{3}\left( 1-3\delta ^{jk}\right).
\end{equation*}
b) The divergence of the $\sigma $\textit{-diagonal vector field} on $M^4$ 
\begin{equation*}
D_{\sigma }=\sigma (x)\frac{\partial }{\partial x^1}+\sigma (x) \frac{%
\partial }{\partial x^2}+\sigma (x)\frac{\partial }{\partial x^3}+ \sigma(x)%
\frac{\partial }{\partial x^4}.
\end{equation*}
has the expression 
\begin{equation*}
\,\mbox{div}\,D_{\sigma }=\sigma _1+\sigma _2+\sigma _3+\sigma _4,
\end{equation*}
where $\sigma _{i}=\partial \sigma /\partial x^{i}$.
\end{lemma}

\noindent Hence, using the general results from \cite{Neagu-Rheon}, we yield:

\begin{proposition}
For the conformal JCM metric (\ref{JCM}), the \textit{energy action
functional} {\small \ 
\begin{eqnarray*}
\mathbb{E}(t,x(t)) &=&\int_{a}^{b}F^2(t,x,y)\sqrt{h_{11}} dt=%
\int_{a}^{b}e^{2\sigma (x)}\cdot G_{11}(y)\cdot h^{11}(t)\sqrt{h_{11}(t)}dt
\\
&=&\int_{a}^{b}e^{2\sigma
(x)}\left(y_1^1y_1^2+y_1^1y_1^3+y_1^1y_1^4+y_1^2y_1^3+
y_1^2y_1^4+y_1^3y_1^4\right) h^{11}\sqrt{h_{11}}dt,
\end{eqnarray*}
}where $y=dx/dt$, produces on the $1$-jet space $J^1(\mathbb{R},M^4)$ the 
\textbf{canonical nonlinear connection} 
\begin{equation}
\Gamma =\left( M_{(1)1}^{(i)}=-\varkappa _{11}^1y_1^{i},\text{ }
N_{(1)j}^{(i)}\right) ,  \label{nlc-B-M}
\end{equation}
where 
\begin{equation*}
N_{(1)j}^{(i)}=\sigma _{j}y_1^{i}+\sigma _{m}y_1^{m}\delta _{j}^{i}+ \left[
\sigma _{i}-\frac{1}{3}\,\mbox{div}\,D_{\sigma }\right]
\left(S_{[1]1}-y_1^{j}\right).
\end{equation*}
\end{proposition}

\begin{myproof}
The Euler-Lagrange equations of the energy action functional $\mathbb{E}$
are 
\begin{equation}
\frac{d^2x^{i}}{dt^2}+2H_{(1)1}^{(i)}\left( t,x^{k},y_1^{k}\right)
+2G_{(1)1}^{(i)}\left( t,x^{k},y_1^{k}\right) =0,\qquad y_1^{k}=\frac{dx^{k}%
}{dt},  \label{Euler-Lagrange=0}
\end{equation}
where we have the local geometrical components 
\begin{equation*}
\left\{%
\begin{array}{lll}
H_{(1)1}^{(i)} & \overset{def}{=} & -\dfrac{1}{2}\varkappa
_{11}^1(t)y_1^{i}\medskip \\ 
G_{(1)1}^{(i)} & \overset{def}{=} & \dfrac{h_{11}g^{ik}}{4}\left[ \dfrac{%
\partial ^2F^2}{\partial x^{m}\partial y_1^{k}}y_1^{m}- \dfrac{\partial F^2}{%
\partial x^{k}}+\dfrac{\partial^2 F^2}{\partial t\partial y_1^{k}}\right.
+\medskip \\ 
&  & \left. +\dfrac{\partial F^2}{\partial y_1^{k}}\varkappa_{11}^1(t)+
2h^{11}\varkappa _{11}^1g_{km}y_1^{m}\right] =\medskip \\ 
& = & \sigma _{m}y_1^{m}y_1^{i}+\left[ \sigma _{i}-\dfrac{1}{3}\,\mbox{div}%
\, D_{\sigma }\right] \cdot G_{11}%
\end{array}%
\right.
\end{equation*}
which determine a \textit{semispray }on the $1$-jet space $J^1(\mathbb{R}%
,M^4)$. Its associated \textit{canonical nonlinear connection} has the
general form \cite{Neagu-Iasi,Neagu-Rheon}) 
\begin{equation*}
\Gamma =\left( M_{(1)1}^{(i)}=2H_{(1)1}^{(i)}=-\varkappa _{11}^1y_1^{i},
\quad N_{(1)j}^{(i)}=\frac{\partial G_{(1)1}^{(i)}}{\partial y_1^{j}}\right).
\end{equation*}
\end{myproof}

\section{Cartan canonical linear connection, d-torsions and d-curvatures}

\hspace{5mm}The canonical nonlinear connection (\ref{nlc-B-M}) is essential
in constructing the dual \textit{adapted bases} of distinguished (d-) vector
fields 
\begin{equation}  \label{a-b-v}
\left\{ \frac{\delta }{\delta t}=\frac{\partial }{\partial t}+
\varkappa_{11}^1y_1^{p}\frac{\partial }{\partial y_1^{p}}\text{ };\text{ } 
\frac{\delta }{\delta x^{i}}=\frac{\partial }{\partial x^{i}}-N_{(1)i}^{(p)} 
\frac{\partial }{\partial y_1^{p}}\text{ };\text{ } \dfrac{\partial }{%
\partial y_1^{i}}\right\} \subset \mathcal{X}(E)
\end{equation}
and distinguished covector fields 
\begin{equation}  \label{a-b-co}
\left\{ dt\text{ };\text{ }dx^{i}\text{ };\text{ }\delta y_1^{i}=dy_1^{i}-
\varkappa_{11}^1y_1^{i}dt+N_{(1)p}^{(i)}dx^{p}\right\} \subset \mathcal{X}%
^{\ast}(E),
\end{equation}
where $E=J^1(\mathbb{R},M^4)$. Note that, under a change of coordinates (\ref%
{tr-rules}), the elements of the adapted bases (\ref{a-b-v}) and (\ref%
{a-b-co}) transform as classical tensors. Consequently, all subsequent
geometrical objects on the $1$-jet space $J^1(\mathbb{R},M^4)$, like Cartan
canonical linear connection, torsion, curvature etc., will be described in
local adapted components.

In this respect, using a general result from \cite{Neagu-Rheon}, by direct
computations, we have the following

\begin{proposition}
The Cartan canonical $\Gamma $-linear connection, produced by the jet
conformal Minkowski metric (\ref{JCM}), has the following adapted local
components: 
\begin{equation}  \label{clc}
C\Gamma =\left( \varkappa _{11}^1,\text{ }G_{j1}^{k}=0,\text{ }L_{jk}^{i}, 
\text{ }C_{j(k)}^{i(1)}=0\right) ,
\end{equation}
where 
\begin{equation*}
L_{jk}^{i}=\delta _{j}^{i}\sigma _{k}+\delta _{k}^{i}\sigma _{j}+
\left(1-\delta _{jk}\right) \sigma _{i}-\frac{1-\delta _{jk}}{3} \,\mbox{div}%
\,D_{\sigma }.
\end{equation*}
\end{proposition}

\begin{myproof}
Using the local derivative operators (\ref{a-b-v}) and the general formulas
which provide the adapted components of the Cartan canonical connection (%
\cite{Neagu-Rheon}), we get 
\begin{equation*}
\left\{%
\begin{array}{l}
G_{j1}^{k}=\frac{g^{km}}{2}\frac{\delta g_{mj}}{\delta t}=0,\quad
L_{jk}^{i}= \frac{g^{im}}{2}\left( \frac{\delta g_{jm}}{\delta x^{k}}+\frac{%
\delta g_{km}}{\delta x^{j}}-\frac{\delta g_{jk}}{\delta x^{m}}%
\right),\smallskip \\ 
C_{j(k)}^{i(1)}=\frac{g^{im}}{2}\left( \frac{\partial g_{jm}}{\partial
y_1^{k}}+\frac{\partial g_{km}}{\partial y_1^{j}}-\frac{\partial g_{jk}}{%
\partial y_1^{m}}\right) =0.%
\end{array}%
\right.
\end{equation*}
\end{myproof}

\begin{remark}
It is straightforward to check the relation $L_{jk}^{i}=\frac{\partial
N_{(1)j}^{(i)}}{\partial y_1^{k}},$ which, considering the homogeneity of
degree $1$ of the local functions $N_{(1)j}^{(i)}$, leads to 
\begin{equation}  \label{equalitie-C}
\frac{\partial N_{(1)j}^{(i)}}{\partial y_1^{m}}y_1^{m}=N_{(1)j}^{(i)}\quad
\Leftrightarrow\quad L_{jm}^{i}y_1^{m}=N_{(1)j}^{(i)}.
\end{equation}
\end{remark}

\begin{proposition}
The Cartan canonical $\Gamma $-linear connection $C\Gamma $ of the jet
conformal Minkowski metric (\ref{JCM}) has a \textbf{single} effective local
torsion d-tensor, namely 
\begin{equation*}
R_{(1)jk}^{(l)}=\mathfrak{R}_{pjk}^{l}y_1^{p},
\end{equation*}
where, using the notations 
\begin{equation*}
\begin{array}{c}
\sigma _{ij}=\dfrac{\partial ^2\sigma }{\partial x^{i}\partial x^{j}}, \quad
\,\mbox{grad}\,\sigma =\left( \sigma _1,\sigma _2,\sigma _3, \sigma_4\right)
,\quad \left\vert \left\vert \,\mbox{grad}\,\sigma \right\vert \right\vert
^2=\sigma _1^2+\sigma _2^2+\sigma _3^2+\sigma_4^2,\smallskip \\ 
\left(\,\mbox{div}\,D_{\sigma }\right) _{i}=\frac{\partial \left( \,%
\mbox{div}\, D_{\sigma }\right) }{\partial x^{i}}=\sigma _{1i}+\sigma _{2i}+
\sigma_{3i}+\sigma _{4i}%
\end{array}%
\end{equation*}
we have 
\begin{eqnarray}
\mathfrak{R}_{ijk}^{l} &=&\delta _{j}^{l}\left( \sigma _{ik}-\sigma_{i}
\sigma _{k}\right) -\delta _{k}^{l}\left( \sigma _{ij}-\sigma _{i}\sigma_{j}
\right) +  \label{FRAKTUR-R} \\
&&+\left( 1-\delta _{ij}\right) \left( \sigma _{lk}-\sigma _{l}\sigma_{k}
\right) -\left( 1-\delta _{ik}\right) \left( \sigma _{lj}-\sigma_{l} \sigma
_{j}\right) +  \notag \\
&&+\frac{1}{3}\left( \,\mbox{div}\,D_{\sigma }\right) \left( \sigma
_{k}-\sigma_{j}+ \delta _{ik}\sigma _{j}-\delta _{ij}\sigma _{k}\right) + 
\notag \\
&&+\left[ \left\vert \left\vert \,\mbox{grad}\,\sigma \right\vert
\right\vert^2- \frac{1}{3}\left( \,\mbox{div}\,D_{\sigma }\right) ^2\right]
\left( \delta _{k}^{l}-\delta _{j}^{l}+\delta _{ik}\delta
_{j}^{l}-\delta_{ij} \delta _{k}^{l}\right) +  \notag \\
&&+\frac{1}{3}\left[ \left( \,\mbox{div}\,D_{\sigma }\right) _{j}-\left( \,%
\mbox{div}\, D_{\sigma }\right) _{k}+\delta _{ij}\left( \,\mbox{div}%
\,D_{\sigma }\right)_{k}- \delta _{ik}\left( \,\mbox{div}\,D_{\sigma
}\right) _{j}\right] {.}  \notag
\end{eqnarray}
\end{proposition}

\begin{myproof}
A general $h$-normal $\Gamma $-linear connection on the 1-jet space $J^1( 
\mathbb{R},M^4)$ is characterized by \textit{eight} effective d-tensors of
torsion (\cite{Neagu-Rheon}). For our Cartan canonical connection \eqref{clc}%
, these reduce only to \textit{one} (the other seven cancel): 
\begin{equation*}
R_{(1)jk}^{(l)}={\dfrac{\delta N_{(1)j}^{(l)}}{\delta x^{k}}}- {\dfrac{%
\delta N_{(1)k}^{(l)}}{\delta x^{j}}}.
\end{equation*}
Using now the expressions of the derivatives $\delta /\delta x^{i}$, formula
(\ref{equalitie-C}) and the $y$-independence $L_{jk}^{i}=L_{jk}^{i}(x)$, we
find 
\begin{equation*}
R_{(1)jk}^{(l)}=\mathfrak{R}_{pjk}^{l}y_1^{p},
\end{equation*}
where 
\begin{equation*}
\mathfrak{R}_{ijk}^{l}:={{\dfrac{\partial L_{ij}^{l}}{\partial x^{k}}}- {%
\dfrac{\partial L_{ik}^{l}}{\partial x^{j}}}%
+L_{ij}^{r}L_{rk}^{l}-L_{ik}^{r}L_{rj}^{l}.}
\end{equation*}
Finally, laborious computations lead to the expression (\ref{FRAKTUR-R}) of
the d-tensor $\mathfrak{R}_{ijk}^{l}$.
\end{myproof}

\begin{proposition}
The Cartan canonical $\Gamma $-linear connection $C\Gamma $ of the jet
conformal Minkowski metric (\ref{JCM}) has a \textbf{single} effective local
curvature d-tensor, namely 
\begin{equation*}
R_{ijk}^{l}=\mathfrak{R}_{ijk}^{l},
\end{equation*}
where $\mathfrak{R}_{ijk}^{l}$ is given by (\ref{FRAKTUR-R}).
\end{proposition}

\begin{myproof}
A general $h$-normal $\Gamma $-linear connection on the 1-jet space $J^1(%
\mathbb{R},M^4)$ is characterized by \textit{five} effective d-tensors of
curvature (\cite{Neagu-Rheon}). For our Cartan canonical connection %
\eqref{clc}, these reduce only to \textit{one} (the other four cancel),
namely 
\begin{equation*}
\begin{array}{lll}
{R_{ijk}^{l}} & \overset{def}{=}\medskip & {{\dfrac{\delta L_{ij}^{l}}{%
\delta x^{k}}}-{\dfrac{\delta L_{ik}^{l}}{\delta x^{j}}}+
L_{ij}^{r}L_{rk}^{l}-L_{ik}^{r}L_{rj}^{l}+C_{i(r)}^{l(1)}R_{(1)jk}^{(r)}=%
\medskip } \\ 
\medskip & = & {{\dfrac{\delta L_{ij}^{l}}{\delta x^{k}}}-{\dfrac{\delta
L_{ik}^{l}}{\delta x^{j}}}+L_{ij}^{r}L_{rk}^{l}-L_{ik}^{r}L_{rj}^{l}}= \\ 
& = & {{\dfrac{\partial L_{ij}^{l}}{\partial x^{k}}}-{\dfrac{\partial
L_{ik}^{l}}{\partial x^{j}}}+L_{ij}^{r}L_{rk}^{l}-L_{ik}^{r}L_{rj}^{l}=} 
\mathfrak{R}_{ijk}^{l}.%
\end{array}%
\end{equation*}
\end{myproof}

\section{Geometrical field model produced by the jet conformal Minkowski
metric}

\subsection{Gravitational-like geometrical model}

\hspace{5mm}From a geometric-physical point of view, on the 1-jet space $J^1(%
\mathbb{R},M^4)$, the jet conformal Minkowski metric (\ref{JCM}) produces
the adapted metrical d-tensor 
\begin{equation}  \label{gravit-pot-B-M}
\mathbb{G}=h_{11}dt\otimes dt+g_{ij}dx^{i}\otimes dx^{j}+h^{11}g_{ij}\delta
y_1^{i}\otimes \delta y_1^{j},
\end{equation}
where $g_{ij}$ is given by (\ref{metric}). This may be regarded as a \textit{%
\textquotedblleft non-isotropic gravitational potential\textquotedblright }.
In such a "physical" context, the nonlinear connection $\Gamma $ (used in
the construction of the distinguished 1-forms $\delta y_1^{i}$) prescribes,
probably, a kind of \textit{\textquotedblleft interaction\textquotedblright }
between $(t)$-, $(x)$- and $(y)$-fields.

We postulate that the non-isotropic gravitational potential $\mathbb{G}$ is
governed by the \textit{geometrical Einstein equations}%
\begin{equation}
\text{Ric }\left( C\Gamma \right) -\frac{\text{Sc }\left( C\Gamma \right) }{2%
}\mathbb{G=}\mathcal{KT},  \label{Einstein-eq-global}
\end{equation}%
where Ric $\left( C\Gamma \right) $ is the \textit{Ricci d-tensor}
associated to the Cartan canonical connection $C\Gamma $ (in Riemannian
sense and using adapted bases), Sc $\left( C\Gamma \right) $ is the \textit{%
scalar curvature}, $\mathcal{K}$ is the \textit{Einstein constant} and $%
\mathcal{T}$ is the intrinsic \textit{stress-energy} d-tensor of matter (%
\cite{MA,BN}).

In this way, working with the adapted basis of vector fields (\ref{a-b-v}),
we can find the local geometrical Einstein equations for the metric (\ref%
{JCM}). Firstly, by direct computations, we find:

\begin{lemma}
The Ricci d-tensor of the Cartan canonical connection $C\Gamma $ of the
metric (\ref{JCM}) has a \textbf{single} effective local Ricci d-tensor,
namely%
\begin{equation}
\begin{array}{lll}
R_{ij} & = & -2\left( \sigma _{ij}-\sigma _{i}\sigma _{j}\right) +\medskip
\\ 
&  & +\dfrac{1-\delta _{ij}}{3}\left[ 3\Delta \sigma +6\left\vert \left\vert
\,\mbox{grad}\,\sigma \right\vert \right\vert ^2-2\left( \,\mbox{div}%
\,D_{\sigma }\right) ^2-\mathfrak{S}\right] ,%
\end{array}
\label{Ricci-local}
\end{equation}
where 
\begin{equation*}
\Delta \sigma =\sigma _{11}+\sigma _{22}+\sigma _{33}+\sigma _{44},\quad 
\mathfrak{S}=\sum_{p,q=1}^4\sigma _{pq}.
\end{equation*}
\end{lemma}

\begin{myproof}
A general $h$-normal $\Gamma $-linear connection on the 1-jet space $J^1( 
\mathbb{R},M^4)$ is characterized by \textit{six} effective Ricci d-tensors (%
\cite{Neagu-Rheon}). For our Cartan canonical connection \eqref{clc}, these
reduce only to \textit{one} (the other five cancel): 
\begin{equation*}
\medskip R_{ij}\overset{def}{=}R_{ijm}^{m}=\mathfrak{R}_{ijm}^{m}.
\end{equation*}
Then, a direct computation gives the expression (\ref{Ricci-local}) of the
Ricci d-tensor $R_{ij}$.
\end{myproof}

\begin{lemma}
The scalar curvature \emph{Sc} $\left( C\Gamma \right) $ of the Cartan
canonical connection $C\Gamma $ of the jet conformal Minkowski metric (\ref%
{JCM}) is given by 
\begin{equation}
R=4e^{-2\sigma }\left[ 3\Delta \sigma +3\left\vert \left\vert \,\mbox{grad}%
\,\sigma \right\vert \right\vert ^2-\left( \,\mbox{div}\,D_{\sigma }\right)
^2-\mathfrak{S}\right] .  \label{Scalar-R}
\end{equation}
\end{lemma}

\begin{myproof}
The general formula for the scalar curvature of the Cartan connection
reduces to (\cite{Neagu-Rheon}) 
\begin{equation*}
\text{Sc }\left( C\Gamma \right) \overset{def}{=}g^{pq}R_{pq}:=R,
\end{equation*}
where $R$ is given by (\ref{Scalar-R}).
\end{myproof}

\smallskip By describing the global geometrical Einstein equations (\ref%
{Einstein-eq-global}) in the adapted basis of vector fields (\ref{a-b-v}),
we find the following important geometrical and physical result (\cite%
{Neagu-Rheon}):

\begin{proposition}
The local \textbf{geometrical Einstein equations} that govern the
non-isotropic gravitational potential $\mathbb{G}$ (produced by the jet
conformal Minkowski metric (\ref{JCM})) are given by: 
\begin{equation}
\medskip R_{ij}-\dfrac{R}{2}g_{ij}=\mathcal{KT}_{ij},  \label{E-1}
\end{equation}
{\small 
\begin{equation}
\left\{%
\begin{array}{llll}
-R\medskip h_{11}=2\mathcal{KT}_{11}, & 0=\mathcal{T}_{1i}, & \medskip 0= 
\mathcal{T}_{i1}, & 0=\mathcal{T}_{(i)1}^{(1)}, \smallskip \\ 
0=\mathcal{T}_{1(i)}^{\text{ }(1)}, & 0=\mathcal{T}_{i(j)}^{\text{ }(1)}, & 
0=\mathcal{T}_{(i)j}^{(1)}, & -Rh^{11}g_{ij}=2\mathcal{KT}_{(i)(j)}^{(1)(1)}.%
\end{array}%
\right.  \label{E-2}
\end{equation}
}
\end{proposition}

\smallskip

\begin{remark}
The Einstein geometrical equations (\ref{E-1}) and (\ref{E-2}) impose that
the stress-energy d-tensor of matter $\mathcal{T}$ be symmetric. In other
words, the stress-energy d-tensor of matter $\mathcal{T}$ must satisfy the
local symmetry conditions 
\begin{equation*}
\mathcal{T}_{AB}=\mathcal{T}_{BA},\quad \forall \text{ }A,B\in \left\{ 1,%
\text{ }i,\text{ }_{(i)}^{(1)}\right\}.
\end{equation*}
Moreover, we must "a priori" have the equality: $\mathcal{T}%
_{(i)(j)}^{(1)(1)} h_{11}=g_{ij}\mathcal{T}_{11}h^{11}$.
\end{remark}

By direct computations, the geometrical Einstein equations (\ref{E-1}) and (%
\ref{E-2}) imply the following identities of the stress-energy d-tensor:%
\footnote{%
Summing over both indices $m,r$ is assumed in \eqref{mr}, over $r$ in %
\eqref{r}, and over $m$ in \eqref{mu} and \eqref{md}.} 
\begin{eqnarray}
&& \mathcal{T}_1^1\overset{def}{=}h^{11}\mathcal{T}_{11}= -\dfrac{R}{2%
\mathcal{K}},\quad \mathcal{T}_1^{m}\overset{def}{=} g^{mr}\mathcal{T}%
_{r1}=0,\quad\mathcal{T}_{(1)1}^{(m)}\overset{def}{=} h_{11}g^{mr}\mathcal{T}%
_{(r)1}^{(1)}=0, \smallskip  \notag \\
&& \mathcal{T}_{i}^1\overset{def}{=}h^{11}\mathcal{T}_{1i}=0,\quad \mathcal{T%
}_{i}^{m}\overset{def}{=}g^{mr}\mathcal{T}_{ri}= \dfrac{1}{\mathcal{K}}%
\left(g^{mr}R_{ri}-\dfrac{R}{2} \mathcal{\delta }_{i}^{m}\right), \smallskip
\notag \\
&& \mathcal{T}_{\text{ \ }(i)}^{1(1)}\overset{def}{=}h^{11} \mathcal{T}%
_{1(i)}^{\text{ }(1)}=0,\quad \mathcal{T}_{(1)(i)}^{(m)(1)}\overset{def}{=}%
h_{11}g^{mr} \mathcal{T}_{(r)(i)}^{(1)(1)}=-\dfrac{R}{2\mathcal{K}} \mathcal{%
\delta }_{i}^{m}, \smallskip  \notag \\
&& \mathcal{T}_{(1)i}^{(m)}\overset{def}{=}h_{11}g^{mr} \mathcal{T}%
_{(r)i}^{(1)}=0,\quad \mathcal{T}_{\text{ \ }(i)}^{m(1)} \overset{def}{=}%
g^{mr}\mathcal{T}_{r(i)}^{\text{ }(1)}=0.  \label{r}
\end{eqnarray}
Consequently, the following local identities for the stress-energy d-tensor
of matter hold good: 
\begin{equation}  \label{mu}
\left\{%
\begin{array}{l}
\mathcal{T}_{1/1}^1+\mathcal{T}_{1|m}^{m}+ \mathcal{T}%
_{(1)1}^{(m)}|_{(m)}^{(1)}=\mathcal{T}_{1/1}^1= -\dfrac{1}{2\mathcal{K}}%
\dfrac{\delta R}{\delta t}\smallskip \\ 
\mathcal{T}_{i/1}^1+\mathcal{T}_{i|m}^{m}+\mathcal{T}%
_{(1)i}^{(m)}|_{(m)}^{(1)}= \mathcal{T}_{i|m}^{m}=\dfrac{1}{\mathcal{K}}%
\left(g^{mr}R_{ri}- \dfrac{R}{2}\mathcal{\delta }_{i}^{m}\right)_{|m}%
\smallskip \\ 
\mathcal{T}_{\text{ \ }(i)/1}^{1(1)}+\mathcal{T}_{\text{ \ }(i)|m}^{m(1)}+ 
\mathcal{T}_{(1)(i)}^{(m)(1)}|_{(m)}^{(1)}= \mathcal{T}%
_{(1)(i)}^{(m)(1)}|_{(m)}^{(1)}= -\dfrac{1}{2\mathcal{K}}\dfrac{\partial R}{%
\partial y_1^{i}},%
\end{array}%
\right.
\end{equation}
where 
\begin{eqnarray}
&& \mathcal{T}_{1/1}^1\overset{def}{=}\dfrac{\delta \mathcal{T}_1^1}{\delta t%
}+ \mathcal{T}_1^1\varkappa _{11}^1-\mathcal{T}_1^1\varkappa _{11}^1= \dfrac{%
\delta \mathcal{T}_1^1}{\delta t},\quad \mathcal{T}_{1|m}^{m}\overset{def}{=}%
\dfrac{\delta \mathcal{T}_1^{m}}{\delta x^{m}}+\mathcal{T}%
_1^{r}L_{rm}^{m}=0, \smallskip  \notag \\
&& \mathcal{T}_{(1)1}^{(m)}|_{(m)}^{(1)}\overset{def}{=}\dfrac{\partial 
\mathcal{T}_{(1)1}^{(m)}}{\partial y_1^{m}}=0,\quad \mathcal{T}_{i/1}^1%
\overset{def}{=}\dfrac{\delta \mathcal{T}_{i}^1}{\delta t}+\mathcal{T}%
_{i}^1\varkappa _{11}^1=0, \smallskip  \notag \\
&& \mathcal{T}_{i|m}^{m}\overset{def}{=}\dfrac{\delta \mathcal{T}_{i}^{m}}{%
\delta x^{m}}+\mathcal{T}_{i}^{r}L_{rm}^{m}- \mathcal{T}_{r}^{m}L_{im}^{r},%
\quad \mathcal{T}_{(1)i}^{(m)}|_{(m)}^{(1)}\overset{def}{=} \dfrac{\partial 
\mathcal{T}_{(1)i}^{(m)}}{\partial y_1^{m}}=0, \smallskip  \notag \\
&& \mathcal{T}_{\text{ \ }(i)/1}^{1(1)}\overset{def}{=}\dfrac{\delta 
\mathcal{T}_{\text{ \ }(i)}^{1(1)}}{\delta t}+ 2\mathcal{T}_{\text{ \ }%
(i)}^{1(1)}\varkappa _{11}^1=0,\quad \mathcal{T}%
_{(1)(i)}^{(m)(1)}|_{(m)}^{(1)}\overset{def}{=}\dfrac{\partial \mathcal{T}%
_{(1)(i)}^{(m)(1)}}{\partial y_1^{m}}, \smallskip  \notag \\
&& \mathcal{T}_{\text{ \ }(i)|m}^{m(1)}\overset{def}{=}\dfrac{\delta 
\mathcal{T}_{\text{ \ }(i)}^{m(1)}}{\delta x^{m}}+ \mathcal{T}_{\text{ \ }%
(i)}^{r(1)}L_{rm}^{m}- \mathcal{T}_{\text{ \ }(r)}^{m(1)}L_{im}^{r}=0.
\label{mr}
\end{eqnarray}
Taking into account that we have the $y$-independence $R=R(x)$, we obtain
the following result:

\begin{corollary}
The stress-energy d-tensor of matter $\mathcal{T}$ must verify the following 
\textbf{conservation geometrical laws}: 
\begin{equation}  \label{md}
\mathcal{T}_{1/1}^1=0,\quad \mathcal{T}_{i|m}^{m}=\dfrac{1}{\mathcal{K}} %
\left[ g^{mr}R_{ri}-\dfrac{R}{2}\mathcal{\delta }_{i}^{m}\right] _{|m},\quad 
\mathcal{T}_{(1)(i)}^{(m)(1)}|_{(m)}^{(1)}=0.
\end{equation}
\end{corollary}

\subsection{Related electromagnetic model considerations}

\hspace{5mm}In the paper \cite{Neagu-Rheon} an electromagnetic geometrical
model was developed, based on a given Lagrangian function $L(t,x,y)$ on the
1-jet space $J^1(\mathbb{R},M^{n})$. In the background of our
electromagnetic geometrical formalism from \cite{Neagu-Rheon}, we work with
an \textit{electromagnetic distinguished $2$-form}\footnote{%
We implicitly assume that the Latin letters run from $1$ to $n$; as well, we
further denote by $\mathcal{A}_{\left\{i,j\right\}}$ -- the alternate sum,
and by $\sum_{\{i,j,k\}}$ -- the cyclic sum.} 
\begin{equation*}
\mathbb{F}=F_{(i)j}^{(1)}\delta y_1^{i}\wedge dx^{j},
\end{equation*}
where 
\begin{equation*}
F_{(i)j}^{(1)}=\frac{h^{11}}{2}\left[
g_{jr}N_{(1)i}^{(r)}-g_{ir}N_{(1)j}^{(r)}+\left(
g_{ir}L_{jm}^{r}-g_{jr}L_{im}^{r}\right) y_1^{m}\right] .
\end{equation*}
The electromagnetic components $F_{(i)j}^{(1)}$ satisfy the following 
\textit{Maxwell geometrical equations} \cite{Neagu-Rheon}: 
\begin{equation*}
\left\{%
\begin{array}{lll}
F_{(i)j/1}^{(1)} & = & \dfrac{1}{2}\,\mathcal{A}_{\left\{ i,j\right\} }
\left\{{\bar D}_{(i)1|j}^{(1)}-D_{(i)m}^{(1)}G_{j1}^{m}+
d_{(i)(m)}^{(1)(1)}R_{(1)1j}^{(m)}-\right. \smallskip \\ 
&  & \left. -\left(C_{j(m)}^{p(1)}R_{(1)1i}^{(m)}-G_{i1|j}^{p}\right)
h^{11}g_{pq}y_1^{q}\right\} ,\smallskip \\ 
\sum_{\{i,j,k\}}F_{(i)j|k}^{(1)} & = & -\dfrac{1}{4}\sum_{\{i,j,k\}} \dfrac{%
\partial^3L}{\partial y_1^{i}\partial y_1^{p}\partial y_1^{m}} \left({\dfrac{%
\delta N_{(1)j}^{(m)}}{\delta x^{k}}}- {\dfrac{\delta N_{(1)k}^{(m)}}{\delta
x^{j}}}\right)y_1^{p},\smallskip \\ 
\sum_{\{i,j,k\}}F_{(i)j}^{(1)}|_{(k)}^{(1)} & = & 0,%
\end{array}%
\right.
\end{equation*}
where 
\begin{equation*}
\begin{array}{ll}
{\bar D}_{(i)1}^{(1)}=\dfrac{h^{11}}{2} \dfrac{\delta g_{im}}{\delta t}%
y_1^{m}, & D_{(i)j}^{(1)}=h^{11}g_{ip}\left(-N_{(1)j}^{(p)}+
L_{jm}^{p}y_1^{m}\right),\smallskip \\ 
\bigskip d_{(i)(j)}^{(1)(1)}=h^{11}\left(g_{ij}+g_{ip}
C_{m(j)}^{p(1)}y_1^{m}\right), & R_{(1)1j}^{(m)}={\dfrac{\delta
M_{(1)1}^{(m)}}{\delta x^{j}}}- {\dfrac{\delta N_{(1)j}^{(m)}}{\delta t},}%
\vspace*{-2mm} \\ 
{\bar D}_{(i)1|j}^{(1)}=\dfrac{\delta {\bar D}_{(i)1}^{(1)}}{\delta x^{j}}-{%
\bar D}_{(m)1}^{(1)}L_{ij}^{m}, & G_{i1|j}^{k}=\dfrac{\delta G_{i1}^{k}}{%
\delta x^{j}}+ G_{i1}^{m}L_{mj}^{k}-G_{m1}^{k}L_{ij}^{m},%
\end{array}%
\end{equation*}
and we have 
\begin{equation*}
\left\{%
\begin{array}{l}
\bigskip F_{(i)j/1}^{(1)}=\dfrac{\delta F_{(i)j}^{(1)}}{\delta t}+
F_{(i)j}^{(1)}%
\varkappa_{11}^1-F_{(m)j}^{(1)}G_{i1}^{m}-F_{(i)m}^{(1)}G_{j1}^{m},\vspace*{%
-2mm} \\ 
F_{(i)j|k}^{(1)}=\dfrac{\delta F_{(i)j}^{(1)}}{\delta x^{k}}-
F_{(m)j}^{(1)}L_{ik}^{m}-F_{(i)m}^{(1)}L_{jk}^{m},\smallskip \\ 
F_{(i)j}^{(1)}|_{(k)}^{(1)}=\dfrac{\partial F_{(i)j}^{(1)}}{\partial y_1^{k}}%
-F_{(m)j}^{(1)}C_{i(k)}^{m(1)}-F_{(i)m}^{(1)}C_{j(k)}^{m(1)}.%
\end{array}%
\right.
\end{equation*}

\begin{example}
The Lagrangian function that governs the movement law of a particle of mass $%
m\neq 0$ and electric charge $e$, which is displaced concomitantly into an
environment endowed both with a gravitational field and an electromagnetic
one, is given by 
\begin{equation}  \label{x05-L-(jet)-ED}
L(t,x^{k},y_1^{k})=mch^{11}(t)\;\varphi _{ij}(x^{k})\;y_1^{i}y_1^{j}+ {\frac{%
2e}{m}}A_{(i)}^{(1)}(t,x^{k})\;y_1^{i}+\mathcal{F}(t,x^{k}),
\end{equation}%
where the semi-Riemannian metric $\varphi _{ij}(x)$ represents the \textit{\
gravitational potential }of the space of events $M$, $A_{(i)}^{(1)}(t,x)$
are the components of a d-tensor on the $1$-jet space $J^1(\mathbb{R},M^{n})$
representing the \textit{electromagnetic potential}\emph{\ }and $\mathcal{F}%
(t,x)$ is a smooth \textit{potential function} on the product manifold $%
\mathbb{R}\times M$. It is important to note that the jet Lagrangian
function (\ref{x05-L-(jet)-ED}) is a natural extension of the Lagrangian
(defined on the tangent bundle) used in electrodynamics by Miron and
Anastasiei \cite{MA}. In our jet Lagrangian formalism applied to ( \ref%
{x05-L-(jet)-ED}), the \textit{electromagnetic components} are given by (%
\cite{Neagu-Rheon}) 
\begin{equation*}
F_{(i)j}^{(1)}=-\frac{e}{2m}\left( \frac{\partial A_{(i)}^{(1)}}{\partial
x^{j}}-\frac{\partial A_{(j)}^{(1)}}{\partial x^{i}}\right) ,
\end{equation*}
and the second set of \textit{Maxwell geometrical equations}\ reduce to the
classical ones \cite{Neagu-Rheon}: 
\begin{equation*}
\sum_{\{i,j,k\}}F_{(i)j|k}^{(1)}=0,
\end{equation*}
where 
\begin{equation*}
F_{(i)j|k}^{(1)}=\frac{\partial F_{(i)j}^{(1)}}{\partial x^{k}}
-F_{(m)j}^{(1)}\gamma _{ik}^{m}-F_{(i)m}^{(1)}\gamma _{jk}^{m}.
\end{equation*}
This fact suggests, in our opinion, applicative credibility for our extended
electromagnetic geometrical model.
\end{example}

On our particular $1$-jet space $J^{1}(\mathbb{R},M^{4})$, the jet conformal
Minkowski metric (\ref{JCM}) (via the Lagrangian function $L=F^{2}$)
produces the electromagnetic $2$-form which due to (\ref{equalitie-C})
trivially vanishes 
\begin{equation*}
\mathbb{F}=0.
\end{equation*}%
In conclusion, the jet conformal Minkowski extended electromagnetic
geometrical model constructed on the $1$-jet space $J^{1}(\mathbb{R},M^{4})$
is trivial. Namely, in our jet geometrical approach, the jet conformal
Minkowski electromagnetism, produced only by the metric \eqref{JCM} alone,
leads to null electromagnetic geometrical components and to tautological
Maxwell-like equations. In our opinion, this fact suggests that the jet
conformal geometrical structure \eqref{JCM} of the $1$-jet space $J^{1}(%
\mathbb{R},M^{4})$ is suitable for modeling gravitation rather than
electromagnetism.{\ }

{\small \noindent\textbf{Authors' addresses}:\newline
\newline
Vladimir B{\scriptsize ALAN}\newline
University Politehnica of Bucharest, Faculty of Applied Sciences,\newline
Department of Mathematics-Informatics I,\newline
Splaiul Independentei 313, RO-060042, Bucharest, Romania.\newline
\textbf{E-mail}: vladimir.balan@upb.ro\newline
\textbf{Website}: http://www.mathem.pub.ro/dept/vbalan.htm\newline
\newline
Mircea N{\scriptsize EAGU}\newline
University Transilvania of Bra\c{s}ov, Faculty of Mathematics and
Informatics,\newline
Department of Algebra, Geometry and Differential Equations,\newline
Blvd. Iuliu Maniu 50, RO-500091, Bra\c{s}ov, Romania.\newline
\textbf{E-mail}: mircea.neagu@unitbv.ro\newline
\textbf{Website}: http://www.2collab.com/user:mirceaneagu }


\begin{thebibliography}{99}
\bibitem{as1} G.S. Asanov, \textit{Finslerian Extension of General Relativity%
}, Reidel, Dordrecht, 1984.

\bibitem{as2} G.S. Asanov, \textit{Jet extension of Finslerian gauge approach%
}, Fortschritte der Physik \textbf{38}, No. \textbf{8} (1990), 571-610.

\bibitem{Atanasiu-Neagu} Gh. Atanasiu, M. Neagu, \textit{On Cartan spaces
with the }$m$\textit{-th root metric }$K(x,p)=\sqrt[m]{a^{i_1i_2\ldots
i_{m}}(x)p_{i_1}p_{i_2}...p_{i_{m}}}$, Hypercomplex Numbers in Geometry and
Physics, No. \textbf{2} (\textbf{12}), Vol. \textbf{6} (2009), 67-73.

\bibitem{BN} V. Balan, M. Neagu, \textit{Jet Single-Time Lagrange Geometry
and Its Applications}, Bucharest, 2010, Preprint.

\bibitem{Bal-Nea-Chernov} V. Balan, M. Neagu, \textit{The }$1$\textit{-jet
generalized Lagrange geometry induced by the rheonomic Chernov metric},
University Politehnica of Bucharest, Scientific Bulletin, Series \textbf{A},
Vol. \textbf{72}, Iss. \textbf{3} (2010), 5-16.

\bibitem{Chernov} V.M. Chernov, \textit{On defining equations for the
elements of associative and commutative algebras and on associated metric
forms}, In: "Space-Time Structure. Algebra and Geometry" (D.G. Pavlov, Gh.
Atanasiu, V. Balan Eds.), Russian Hypercomplex Society, Lilia Print, Moscow,
2007, 189-209.

\bibitem{Garasko-Book} G.I. Garas'ko, \textit{Foundations of Finsler
Geometry for Physicists}, Tetru Eds, Moscow, 2009 (in Russian).

\bibitem{Garasko-Pavlov} G.I. Garas'ko, D.G. Pavlov, \textit{The notions of
distance and velocity modulus in the linear Finsler spaces}, In: "Space-Time
Structure. Algebra and Geometry" (D.G. Pavlov, Gh. Atanasiu, V. Balan Eds.),
Russian Hypercomplex Society, Lilia Print, Moscow, 2007, 104-117.

\bibitem{Mikhailov} R.V. Mikhailov, \textit{On some questions of four
dimensional topology: a survey of modern research}, Hypercomplex Numbers in
Geometry and Physics, No. \textbf{1} (\textbf{1}), Vol. \textbf{1} (2004),
100-103.

\bibitem{MA} R. Miron, M. Anastasiei, \textit{The Geometry of Lagrange
Spaces: Theory and Applications}, Kluwer Academic Publishers, 1994.

\bibitem{Neagu-Iasi} M. Neagu, \textit{Jet geometrical objects depending on
a relativistic time}, An. \c{S}t. Univ. "Al.I. Cuza" Ia\c{s}i (S.N.). Mat.,
Tomul \textbf{LVI}, f. \textbf{2} (2010), 407-428.

\bibitem{Nc} M. Neagu, \textit{Riemann-Lagrange Geometry on }$1$ \textit{%
-Jet Spaces}, Matrix Rom, Bucharest, 2005.

\bibitem{Neagu-Rheon} M. Neagu, \textit{The geometry of relativistic
rheonomic Lagrange spaces}, BSG Proceedings \textbf{5}, Geometry Balkan
Press, Bucharest (2001), 142-168.

\bibitem{Olver} P.J. Olver, \textit{Applications of Lie Groups to
Differential Equations}, Springer-Verlag, 1986.

\bibitem{Pavlov} D.G. Pavlov, \textit{Four-dimensional time}, Hypercomplex
Numbers in Geometry and Physics, No. \textbf{1} (\textbf{1}), Vol. \textbf{1}
(2004), 31-39.

\bibitem{pa2} D.G. Pavlov, \emph{Generalization of scalar product axioms},
Hypercomplex Numbers in Geometry and Physics, No. \textbf{1} (\textbf{1}),
Vol. \textbf{1} (2004), 5--18.

\bibitem{pa3} D.G. Pavlov, S.S. Kokarev, \emph{Conformal gauges of the
Berwald-Moor Geometry and their induced non-linear symmetries} (in Russian),
Hypercomplex Numbers in Geometry and Physics No. \textbf{2} (\textbf{10}),
Vol. \textbf{5} (2008), 3--14.

\bibitem{Shimada} H. Shimada, \emph{On Finsler spaces with the metric\newline
$L(x,y)=\sqrt[m]{a_{i_1i_2\ldots i_m}(x)y^{i_1}y^{i_2}...y^{i_m}}$}, Tensor
N. S. \textbf{33} (1979), 365-372.
\end{thebibliography}
\end{document}